\documentclass[article,12pt]{amsart}
\usepackage{amsfonts}
\usepackage{amsmath}
\usepackage{graphicx}
\usepackage{color}
\usepackage{epsfig}

\numberwithin{equation}{section}
\theoremstyle{plain}
\newtheorem{thm}{Theorem}[section]
\theoremstyle{remark}
\newtheorem{rem}{Remark}[section]

\newcommand{\cA}{\mathcal{A}}
\newcommand{\cB}{\mathcal{B}}
\newcommand{\cC}{\mathcal{C}}

\newcommand{\cO}{\mathcal{O}}
\newcommand{\cP}{\mathcal{P}}
\newcommand{\cQ}{\mathcal{Q}}

\newcommand{\ind}{\mbox{1}\kern-.25em \mbox{I}}
\font\calcal=cmsy10 scaled\magstep1
\def\build#1_#2^#3{\mathrel{\mathop{\kern 0pt#1}\limits_{#2}^{#3}}}

\def\videbox{\mathbin{\vbox{\hrule\hbox{\vrule height1ex \kern.5em
\vrule height1ex}\hrule}}}
\newtheorem{lem}{Lemma}[section]

\email{Bernard.Bercu@math.u-bordeaux1.fr}
\email{Victor.Vazquez@siu.buap.mx}

\begin{document}
\title[On the usefulness of persistent excitation in ARX adaptive tracking]
{On the usefulness of persistent excitation in ARX adaptive tracking}
\author{Bernard Bercu}
	\address{ Universit\'e Bordeaux 1, Institut de Math\'ematiques de Bordeaux,
	UMR 5251, 351 cours de la lib\'eration, 33405 Talence cedex, France.}
\author{Victor Vazquez}
\address{ Universidad Aut\'onoma de Puebla, Facultad de Ciencias Fisico
Matem\'aticas, Avenida San Claudio y Rio Verde, 72570 Puebla, Mexico.}
\thanks{This work has been supported by INRIA, by CONACYT, 
and by the ECOS Scientific Cooperation Programme.}

\begin{abstract}
The usefulness of persistent excitation is well-known in the
control community. Thanks to a persistently excited adaptive tracking control,
we show that it is possible to avoid the strong controllability assumption
recently proposed in the multidimensional ARX framework.
We establish the almost sure convergence for both least squares
and weighted least squares estimators of the unknown parameters. 
A central limit theorem and a law of iterated logarithm are also provided.
All this asymptotical analysis is related to the Schur complement of a suitable limiting
matrix.
\end{abstract}

\maketitle


\section{INTRODUCTION}


The concept of persistent excitation is well-known in the
control community. Since the pioneers works of Anderson
\cite{Anderson} and Moore \cite{Moore}, this concept has been successfullly used
in a large variety of fields of application
going from economics \cite{Aggelogiannaki}, \cite{Corrado}, 
to adaptive or learning control \cite{Chengyu}, \cite{Farrell},
\cite{Levanoy1}, \cite{Levanoy2}, or mechanical engineering 
and robotics \cite{Abouelsoud}, \cite{Fang}, and \cite{Huang}.
In this paper, we use a persistently excited adaptive tracking control
in the multidimensional ARX framework. It allows us to avoid the strong controllability 
assumption recently proposed by Bercu and Vazquez \cite{BercuVazq}, \cite{BercuVazquez}.  
More precisely, we shall establish the almost sure convergence
for both least squares (LS) and weighted least squares (WLS) estimators of 
the unknown parameters of ARX model. The asymptotic normality as well as a law of iterated logarithm are also provided.
Consider the $d$-dimensional autoregressive process with adaptive control of
order $(p,q)$, $\mbox{ARX}_{d}(p,q)$ for short, given for all $n\geq 0$ by 
\begin{equation}  
\label{ARX}
A(R)X_{n+1}=B(R)U_{n}+\varepsilon_{n+1}
\end{equation}
where $R$ stands for the shift-back operator and $X_{n},U_{n}$ and 
$\varepsilon _{n}$ are the system output, input and driven noise,
respectively. The polynomials $A$ and $B$ are given for all $z\in \mathbb{C}$
by 
\begin{eqnarray*}
A(z) &=&I_{d}-A_{1}z-\cdots -A_{p}z^{p}, \\
B(z) &=&I_{d}+B_{1}z+\cdots +B_{q}z^{q},
\end{eqnarray*}
where $A_{i}$ and $B_{j}$ are unknown square matrices of order $d$ and 
$I_{d} $ is the identity matrix. Relation (\ref{ARX}) may be rewritten 
in the compact form
\begin{equation}  
\label{MOD}
X_{n+1}=\theta ^{t}\Phi _{n}+U_{n}+\varepsilon_{n+1}
\end{equation}
where the regression vector $\Phi _{n}=\left( X_{n}^{p},U_{n-1}^{q}\right)^{t}$ 
with 
\begin{eqnarray*}
X_{n}^{p} &=&(X_{n}^{t},\ldots ,X_{n-p+1}^{t}), \\
U_{n}^{q} &=&(U_{n}^{t},\ldots ,U_{n-q+1}^{t}),
\end{eqnarray*}
and the unknown parameter $\theta$ is given by 
\begin{equation*}
\theta ^{t}=(A_{1},\ldots ,A_{p},B_{1},\ldots ,B_{q}).
\end{equation*}
In all the sequel, we shall assume that the driven noise 
$(\varepsilon _{n})$ is a martingale difference sequence adapted
to the filtration $\mathbb{F}=(\mathcal{F}_{n})$
where $\mathcal{F}_{n}$
stands for the $\sigma $-algebra of the events occurring up to time $n$. Moreover, we
also assume that, for all $n\geq 0$, 
$\mathbb{E}[\varepsilon_{n+1}\varepsilon _{n+1}^{t}|\mathcal{F}_{n}]=\Gamma $ a.s. where $\Gamma $
is a positive definite deterministic covariance matrix. In addition, we suppose that the driven noise 
$(\varepsilon _{n})$ satisfies the strong law of large numbers i.e. if 
\begin{equation}
\Gamma_{n}=\frac{1}{n}\sum_{k=1}^{n}\varepsilon_{k}\varepsilon_{k}^{t},
\end{equation}
then the sequence $(\Gamma_{n})$ converges to $\Gamma$ a.s. That is the case if, for
example, $(\varepsilon _{n})$ is a white noise or if $(\varepsilon _{n})$ has a
finite conditional moment of order $>2$. \\ \par
The paper is organized as follows. Section $\!2$ deals with the parameter
estimation and the persistently excited adaptive tracking control. Section $\!3$ is devoted to the
introduction of the Schur complement approach together with some
linear algebra calculations. In Section $\!4$, we propose
some usefull almost sure convergence properties together with
a central limit theorem (CLT) and a law of iterated logarithm (LIL) for both LS
and WLS estimators. Some numerical simulations are also provided in Section $\!5$.
Finally, a short conclusion is given in Section $\! 6$.


\section{Estimation and Adaptive control}


In the ARX tracking framework, we must deal with two objectives simultaneously.
On the one hand, it is necessary to estimate the unknown parameter $\theta$.
On the other hand, the output $(X_n)$ has to track, step by step, a predictable 
reference trajectory $(x_n)$.
First, we focus our attention on the estimation of the parameter 
$\theta$. We shall make use of the WLS algorithm
which satisfies, for all $n\geq 0$, 
\begin{equation}  
\label{WLS}
\widehat{\theta}_{n+1}=\widehat{\theta}_{n}+a_nS_{n}^{-1}(a)\Phi_{n}
\left(X_{n+1}-U_{n}-\widehat{\theta}_{n}^{\,t}\Phi_{n}\right){\! }^{t}
\end{equation}
where the initial value $\hat{\theta}_{0}$ may be arbitrarily chosen and 
\begin{equation*}
S_{n}(a)=\sum_{k=0}^{n} a_k\Phi_{k}\Phi_{k}^{t}+I_\delta
\end{equation*}
where the identity matrix $I_\delta$ with $\delta=d(p+q)$ is added in order
to avoid the useless invertibility assumption. The choice of the weighted
sequence $(a_{n})$ is crucial. If 
\begin{equation*}
a_{n}=1
\end{equation*}
we find the standard LS estimator, while if $\gamma \!> \! 0$, 
\begin{equation*}
a_{n}=\Bigl(\frac{1}{\log s_{n}} \Bigr)^{1+\gamma} \hspace{0.5cm}\text{with}%
\hspace{0.5cm} s_n=\sum_{k=0}^n\parallel\Phi_k\parallel^2,
\end{equation*}
we obtain the WLS estimator introduced by Bercu and Duflo 
\cite{BercuDuflo}, \cite{Bercu1}.
Next, we are concern with the choice of the adaptive control sequence $(U_n)$. The
crucial role played by $U_n$ is to regulate the dynamic of the process 
$(X_n) $ by forcing $X_n$ to track a predictable reference
trajectory $(x_n)$. We propose to make use of the persistently excited adaptive 
tracking control given, for all $n \geq 0$, by 
\begin{equation}  \label{CONTROL}
U_n = x_{n+1}-\widehat{\theta}_n^{\,t}\,\Phi_n+\xi_{n+1}
\end{equation}
where $(\xi_n)$ is an exogenous noise of dimension $d$, adapted to $\mathbb{F}$, with mean $0$ and positive definite covariance matrix $\Delta$. In addition, we assume that $(\xi_n)$ is independent of  $(\varepsilon_n)$, of $(x_n)$, and of the initial state of the system. Moreover, we suppose 
that $(\xi_n)$ satisfies the strong law of large numbers. Consequently, if  
\begin{equation}  \label{DELTAN}
\Delta_n = \frac{1}{n}\sum_{k=1}^{n}(\varepsilon_{k}+\xi_k)(\varepsilon_{k}+\xi_k)^{t},
\end{equation}
then the sequence $(\Delta_n)$ converges to $\Gamma+\Delta$ a.s. 
By substituting (\ref{CONTROL}) into (\ref{MOD}), we obtain the closed-loop
system 
\begin{equation}  \label{CLS}
X_{n+1} - x_{n+1}= \pi_n + \varepsilon_{n+1}+\xi_{n+1}
\end{equation}
where the prediction error $\pi_n = (\theta - \widehat\theta_n)^{\,t}\Phi_n$.
Furthermore, we assume in all the sequel that the reference trajectory $(x_n)$
satisfies 
\begin{equation}  
\label{CT}
\sum_{k=1}^{n} \parallel x_{k} \parallel^{2} =o(n) \hspace{1cm} \text{a.s.}
\end{equation}
Finally, let $(C_{n})$ be the
average cost matrix sequence defined by 
\begin{equation*}
C_{n}=\frac{1}{n}\sum_{k=1}^{n}(X_{k}-x_{k})(X_{k}-x_{k})^{t}.
\end{equation*}
The tracking is said to be residually optimal if $(C_{n})$ converges to $\Gamma+\Delta$
a.s.


\section{On the Schur Complement}


In all the sequel, we shall make use of the well-known causality assumption
on $B$. More precisely, we assume that for all $z\in \mathbb{C}$ with 
$|z|\leq 1$ 
\begin{equation}
\det(B(z))\neq 0. 
\end{equation}
In other words, the polynomial $\det(B(z))$ only has zeros with modulus $> 1$. 
Consequently, if $r>1$ is strictly less than the smallest modulus of the
zeros of $\det(B(z))$, then $B(z)$ is invertible in the ball with center
zero and radius $r$ and $B^{-1}(z)$ is a holomorphic function (see e.g. \cite {Duflo} page 155).
Hence, for all $z\in \mathbb{C}$ with $|z|\leq r$, we have 
\begin{equation}  \label{INVB}
B^{-1}(z)=\sum_{k=0}^\infty D_kz^k.
\end{equation}
where all the matrices $D_k$ can be explicitly calculated via the recursive equations
$D_0=I_d$ and, for all $k\geq 1$
\begin{eqnarray}
D_{k} &=&-\sum_{j=0}^{k-1}D_{j}B_{k-j}\quad \text{if}  \quad k\leq q,  \label{inv1} \\
D_{k} &=&-\sum_{j=1}^{q}D_{k-j}B_{j} \quad \text{if}  \quad k>q. \label{inv2}
\end{eqnarray}
In a similar way, for all $z\in \mathbb{C}$ such that $|z|\leq r$, we shall
denote 
\begin{equation}  \label{DEFP}
P(z)=B^{-1}(z)(A(z)-I_d)=\sum_{k=1}^\infty P_kz^k.
\end{equation}
All the matrices $P_k$ may be explicitly calculated as functions of the
matrices $A_{i}$ and $B_{j}$. As a matter of fact, for all $k\geq 1$
\begin{eqnarray}
P_{k} &=&-\sum_{j=0}^{k-1}D_{j}A_{k-j} \quad \text{if}  \quad k\leq p,  \label{P1} \\
P_{k} &=&-\sum_{j=1}^{p}D_{k-j}A_{j} \quad \text{if}  \quad k>p.   \label{P2}
\end{eqnarray}
For all $1\leq i \leq q$, denote by $H_i$ be the square matrix of order $d$ 
\begin{equation*}
H_i=\sum_{k=i}^\infty P_k\Gamma P_{k-i+1}^{t}+\sum_{k=i-1}^\infty Q_k\Delta Q_{k-i+1}^{t}.
\end{equation*}
where,  for all $k \geq 0$, $Q_k=D_k+P_k$ with $Q_0=I_d$. In addition, let $H$ be the symmetric square matrix of order $dq$ 
\begin{equation}  \label{DEFH}
H=\left( 
\begin{array}{ccccc}
H_1 & H_2 & \cdots & H_{q-1} & H_q \\ 
H_2^t & H_1 & H_2 & \cdots & H_{q-1} \\ 
\cdots & \cdots & \cdots & \cdots & \cdots \\ 
H_{q-1}^t & \cdots & H_2^t & H_1 & H_2 \\ 
H_q^t & H_{q-1}^t & \cdots & H_2^t & H_1%
\end{array}
\right).
\end{equation}
For all $0\leq i \leq p-1$, let $K_i=P_i\Gamma+Q_i\Delta $ with $K_0=\Delta$ 
and denote by $K$ the rectangular matrix of dimension $dq\times dp$ given, if $p\geq q$, by 
\begin{equation*}
K=\left( 
\begin{array}{ccccccc}
K_0 & K_1 & K_2 & \cdots & \cdots & K_{p-2} & K_{p-1} \\ 
0 & K_0 & K_1 & \cdots & \cdots & K_{p-3} & K_{p-2} \\ 
\cdots & \cdots & \cdots & \cdots & \cdots & \cdots & \cdots \\ 
0 & \cdots & K_0 & K_1 & K_2 & \cdots & K_{p-q+1} \\ 
0 & \cdots & \cdots & K_0 & K_1 & \cdots & K_{p-q}%
\end{array}
\right)
\end{equation*}
while, if $p\leq q$, by 
\begin{equation*}
K=\left( 
\begin{array}{ccccc}
K_0 & K_1 & \cdots & K_{p-2} & K_{p-1} \\ 
0 & K_0 & K_1 & \cdots & K_{p-2} \\ 
\cdots & \cdots & \cdots & \cdots & \cdots \\ 
0 & \cdots & 0 & K_0 & K_{1} \\ 
0 & 0 & \cdots & 0 & 0 \\ 
\cdots & \cdots & \cdots & \cdots & \cdots \\ 
0 & 0 & \cdots & 0 & 0%
\end{array}
\right).
\end{equation*}
Finally, let $L$ be the block diagonal matrix of order $dp$ 
\begin{equation}  \label{DEFL}
L=\left( 
\begin{array}{ccccc}
\Gamma+\Delta & 0 & \cdots & 0 & 0 \\ 
0 & \Gamma+\Delta & 0 & \cdots & 0 \\ 
\cdots & \cdots & \cdots & \cdots & \cdots \\ 
0 & \cdots & 0 & \Gamma+\Delta & 0 \\ 
0 & 0 & \cdots & 0 & \Gamma+\Delta
\end{array}
\right).
\end{equation}
Denote by $\Lambda$ the symmetric square matrix of order $\delta$ 
\begin{equation}  \label{DEFLAMBDA}
\Lambda=\left( 
\begin{array}{cc}
L & K^t \\ 
K & H%
\end{array}
\right).
\end{equation}
This lemma is the keystone of all our asymptotic results. \\
\begin{lem}
\label{MAINLEMMA} Let $S$ be the Schur complement of $L$ in $\Lambda$ 
\begin{equation}  \label{SCHUR}
S=H-KL^{-1}K^t.
\end{equation}
If $B$ is causal, then $S$ and $\Lambda$ are invertible
and 
\begin{equation}  \label{INVLAMBDA}
\Lambda^{-1}\!=\! \left( 
\begin{array}{cc}
\!\!L^{-1} \!+\! L^{-1}K^tS^{-1}KL^{-1} & \!- L^{-1}K^tS^{-1}\!\! \\ 
\!- S^{-1}KL^{-1} & S^{-1}
\end{array}
\right).
\end{equation}
\end{lem}
\ \\
\begin{proof}
The proof is given in Appendix\,A.
\end{proof}
\begin{rem}
One can see the usefulness of persistent excitation in ARX tracking.
As we make use of a persistently excited adaptive 
tracking control given, it is possible to get ride of the strong controllability 
assumption recently proposed by Bercu and Vazquez \cite{BercuVazq}, \cite{BercuVazquez}.
On the other hand, we will see in the next section that the tracking is not optimal
but it is residually optimal. It is necessary to make a compromise
between estimation and tracking optimality.
\end{rem}


\section{MAIN RESULTS}


Our first result concerns to the a.s. asymptotic properties of the LS estimator.
\begin{thm}
\label{ASPLS} Assume that $B$ is causal and that $(\varepsilon_n)$ has finite 
conditional moment of order $>2$. Then, for the LS estimator, we have 
\begin{equation}  \label{TH11}
\lim_{n\rightarrow \infty} \frac{S_{n}}{n} = \Lambda \hspace{1cm} 
\text{a.s.}
\end{equation}
where the limiting matrix $\Lambda$ is given by (\ref{DEFLAMBDA}). In
addition, the tracking is residually optimal 
\begin{equation}  \label{TH12}
\parallel C_{n}-\Delta_{n} \parallel = \mathcal{O} \left( \frac{\log n}{n}
\right) \hspace{1cm}\text{a.s.}
\end{equation}
Finally, $\widehat{\theta}_{n}$ converges almost surely to $\theta$ 
\begin{equation}  \label{TH14}
\parallel \widehat{\theta}_{n}-\theta \parallel^{2}= \mathcal{O} \left(\frac{%
\log n}{n} \right) \hspace{1cm}\text{a.s.}
\end{equation}
\end{thm}
\begin{proof}
The proof is given in Appendix\,B.
\end{proof}
\vspace{2ex}
Our second result is related to the almost sure properties of the WLS
estimator.
\begin{thm}
\label{ASPWLS} 
Assume that $B$ is causal. In addition, suppose that either $(\varepsilon_n)$ is a white
noise or $(\varepsilon_n)$ has finite conditional moment of order $>2$.
Then, for the WLS estimator, we have 
\begin{equation}  \label{TH21}
\lim_{n\rightarrow \infty} (\log n)^{1+\gamma} \frac{S_{n}(a)}{n} = \Lambda 
\hspace{1cm}\text{a.s.}
\end{equation}
where the limiting matrix $\Lambda$ is given by (\ref{DEFLAMBDA}). In
addition, the tracking is residually optimal 
\begin{equation}  \label{TH22}
\parallel C_{n}-\Delta_{n} \parallel= o\left( \frac{(\log n)^{1+\gamma}}{n}
\right) \hspace{1cm}\text{a.s.}
\end{equation}
Finally, $\widehat{\theta}_{n}$ converges almost surely to $\theta$ 
\begin{equation}  \label{TH23}
\parallel \widehat{\theta}_{n}-\theta \parallel^{2}= \mathcal{O} \left( 
\frac{(\log n)^{1+\gamma}}{n} \right) \hspace{1cm}\text{a.s.}
\end{equation}
\end{thm}
\begin{proof}
The proof is given in Appendix\,C.
\end{proof}
\vspace{2ex}
Finally, we present the CLT and the LIL for both LS and WLS estimators.
\begin{thm}
\label{CLTLIL} 
Assume that $B$ is causal and that $(\varepsilon_n)$ and $(\xi_n)$ 
have both finite conditional moments of
order $\alpha>2$. In addition, suppose that $(x_n)$ satisfies 
for some $2< \beta < \alpha$ 
\begin{equation}  \label{REGNORM}
\sum_{k=1}^{n}\parallel x_{k}\parallel^{\beta}=\mathcal{O}(n) \hspace{0.5cm} 
\text{a.s.}
\end{equation}
Then, the LS and WLS estimators share the same central limit theorem 
\begin{equation}  \label{TH31}
\sqrt{n} (\widehat{\theta}_{n}-\theta ) \build{\longrightarrow}_{}^{{%
\mbox{\calcal L}}} \mathcal{N}(0,\Lambda^{-1}\otimes\Gamma)
\end{equation}
where the inverse matrix $\Lambda^{-1}$ is given by (\ref{INVLAMBDA}) and
the symbol $\otimes$ stands for the matrix Kronecker product. In addition,
for any vectors $u \in \mathbb{R}^{d}$ and $v \in \mathbb{R}^{\delta}$, they
also share the same law of iterated logarithm 
\begin{eqnarray}  \label{TH32}
\limsup_{n \rightarrow \infty} \left(\frac{n}{2 \log\log n} \right)^{1/2}
\!\!v^{t}(\widehat{\theta}_{n}-\theta)u 
&=& - \liminf_{n \rightarrow \infty}
\left(\frac{n}{2 \log\log n}\right)^{1/2} \!\!v^{t}(\widehat{\theta}
_{n}-\theta)u  \notag \\
&=& 
\Bigl(v^{t}\Lambda^{-1}v \Bigr)^{1/2}\Bigl(u^{t}\Gamma u \Bigr)^{1/2} 
\hspace{0.5cm}\text{a.s.}
\end{eqnarray}
In particular, 
\begin{equation*}  
\left(\frac{\lambda_{min}\Gamma} {\lambda_{max}\Lambda} \right)\leq 
\limsup_{n \rightarrow \infty} \ \left(\frac{n}{2 \log\log n} \right)
\parallel \hat{\theta}_{n}-\theta \parallel^{2}  \leq 
 \left(\frac{\lambda_{max}\Gamma} {\lambda_{min}\Lambda}\right) 
 \hspace{0.2cm}\text{a.s.}
\end{equation*}
where $\lambda_{min} \Gamma$ and $\lambda_{max} \Gamma$ are the minimum and
the maximum eigenvalues of $\Gamma$.
\end{thm}
\begin{proof}
The proof is given in Appendix\,D.
\end{proof}


\section{NUMERICAL SIMULATIONS}


The goal of this section is to illustrate via some numerical experiments the main results of this paper. 
In order to keep this section brief, we consider a causal $ARX_d(p,q)$ model 
in dimension $d=2$ with $p=1$ and $q=1$. Moreover, the reference trajectory $(x_n)$ is chosen to be
identically zero and the driven and exogenous noises $(\varepsilon_n)$ and $(\xi_n)$ are Gaussian 
$\mathcal{N}(0,1)$ white noises. Finally our numerical simulations are based on $M=500$ realizations of sample size $N=1000$. Consider the 
$ARX_{2}(1,1)$ model 
\begin{equation*}
X_{n+1}=A X_n+U_{n}+B U_{n-1} + \varepsilon _{n+1}
\end{equation*}
where 
\begin{equation*}
A =\left( 
\begin{array}{cc}
2 & 0 \\ 
0 & 0
\end{array}
\right) \hspace{0.5cm} \text{and} \hspace{0.5cm} 
B =\frac{1}{4}\left( 
\begin{array}{cc}
3 & \,0 \\ 
0 & \!-2
\end{array}
\right).
\end{equation*}
First of all, it is easy to see that this $ARX_{2}(1,1)$ process is
not strongly controllable \cite{BercuVazq}, \cite{BercuVazquez}, because $\det(A)=0$. Consequently, 
if we use an adaptive tracking control $(U_n)$ without
persistent excitation $(\xi_n)$, then only the matrix $A$ and 
the first diagonal term of the matrix $B$
can be properly estimated as one can see in Figure 1. 
\begin{figure}[htp]
\vspace{-5cm}
\centering
\includegraphics[width=2.7in, height=8in]{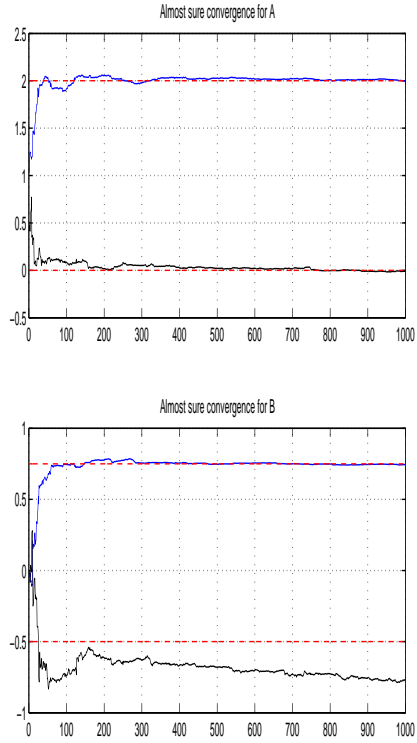}
\vspace{-5.5cm}
\caption{Almost Sure Convergence}
\end{figure}
\ \\
Next, we make use of the persistently excited adaptive tracking control
given by
\begin{equation*}  
U_n = -\widehat{\theta}_n^{\,t}\,\Phi_n+\xi_{n+1}.
\end{equation*}
For all $k\geq 1$, we have $D_k=(-B)^{k}$ and $P_k=-(-B)^{k-1}A$
which clearly implies that 
$$Q_k=-(-B)^{k-1}(A+B).$$
Since the matrices $A$ and $B$ are both diagonal, we find that
\begin{eqnarray*}
H&=&\sum_{k=1}^\infty P_k^2+ \sum_{k=0}^\infty Q_k^2,\\
&=&I_2+\!\sum_{k=1}^\infty B^{k-1}A^2B^{k-1}+\!\sum_{k=1}^\infty B^{k-1}(A+B)^2B^{k-1},\\
&=&I_2+(A^2+(A+B)^2)\sum_{k=0}^\infty B^{2k},\\
&=&I_2+(A^2+(A+B)^2)(I_2 - B^2)^{-1}.
\end{eqnarray*}
Consequently, we obtain that
\begin{equation*}
H=\frac{1}{21}\left( 
\begin{array}{cc}
576 & 0 \\ 
0 & 28
\end{array}
\right).
\end{equation*}
Therefore, the limiting matrix $\Lambda$ given 
by (\ref{DEFLAMBDA}) is
\begin{equation*}
\Lambda=\frac{1}{21}\left( 
\begin{array}{cccc}
42 & 0 & 21 & 0 \\ 
0 & 42 & 0 & 21 \\
21 & 0 & 576 & 0 \\
0 & 21 & 0 & 28
\end{array}
\right).
\end{equation*}
It is not hard to see that
$\det(\Lambda)=89.7619$.
One can observe in Figure 2 the almost sure convergence of the LS estimator 
$\widehat{\theta}_n$ to the four diagonal coordinates of $\theta$. One can conclude
that $\widehat{\theta}_n$ performs very well in the estimation of $\theta$. 
\ \\
\begin{figure}[htp]
\vspace{-0.5cm}
\centering
\includegraphics[width=2.2in, height=4in]{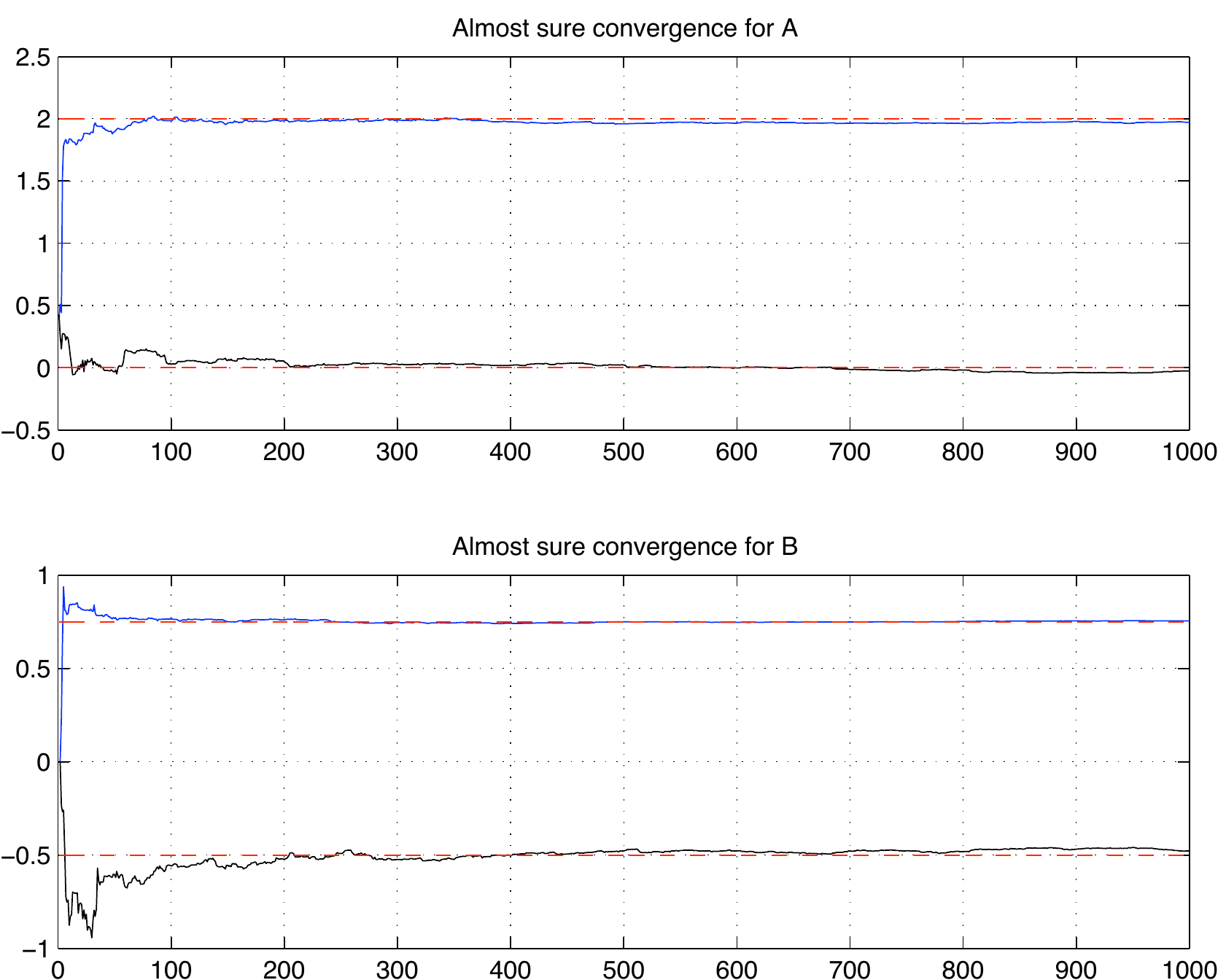}
\caption{Almost Sure Convergence}
\end{figure}
\ \\
Figure 3 shows the CLT for the four coordinates of 
\begin{equation*}
Z_N=\sqrt{N}\Lambda^{1/2}(\widehat{\theta}_N - \theta).
\end{equation*}
One can realize that each component of $Z_N$ has $\mathcal{N}(0,1)$
distribution as expected.
\begin{figure}[htp]
\centering
\includegraphics[width=2.2in, height=4in]{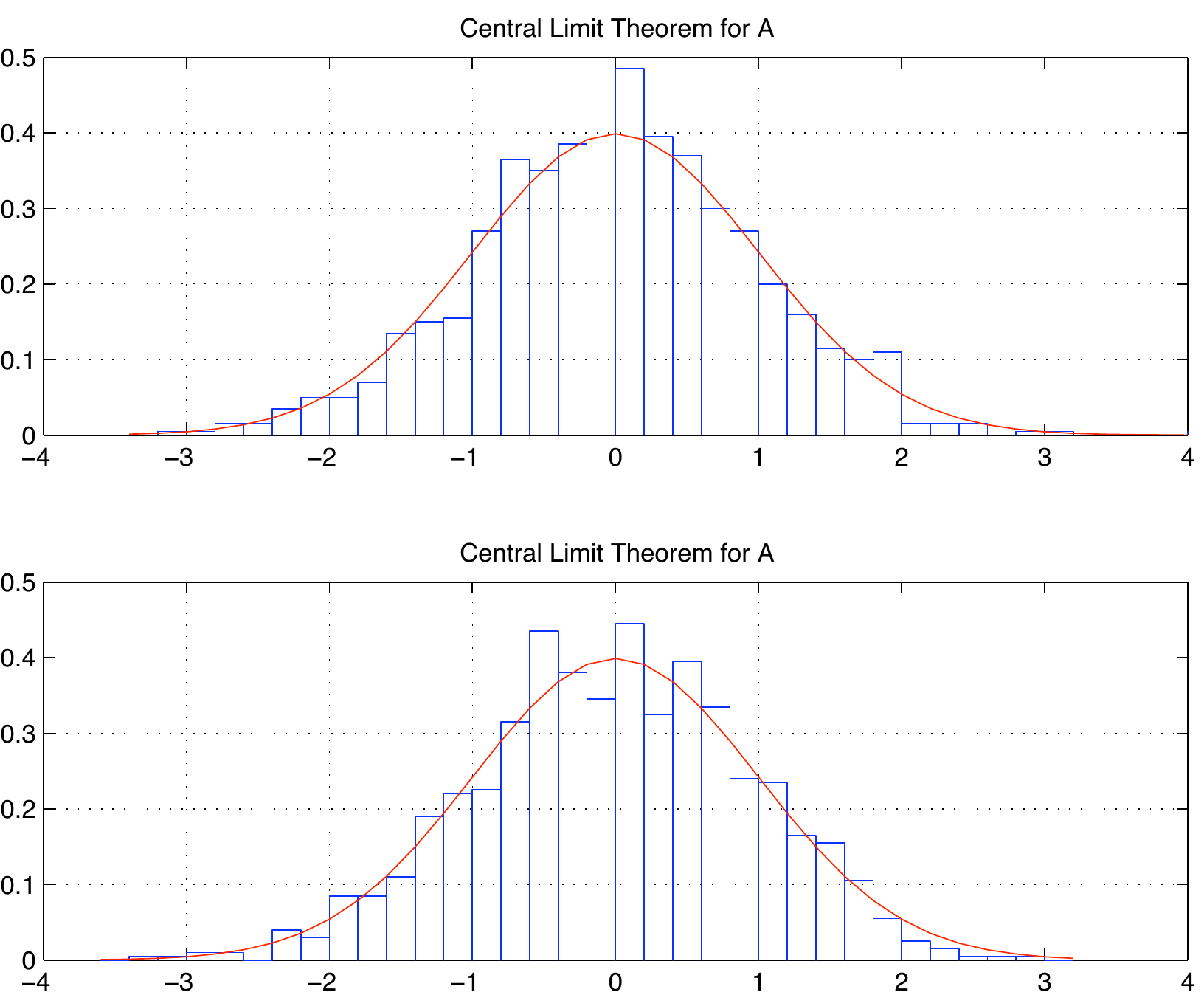}
\end{figure}
\begin{figure}[htp]
\centering
\includegraphics[width=2.2in, height=4in]{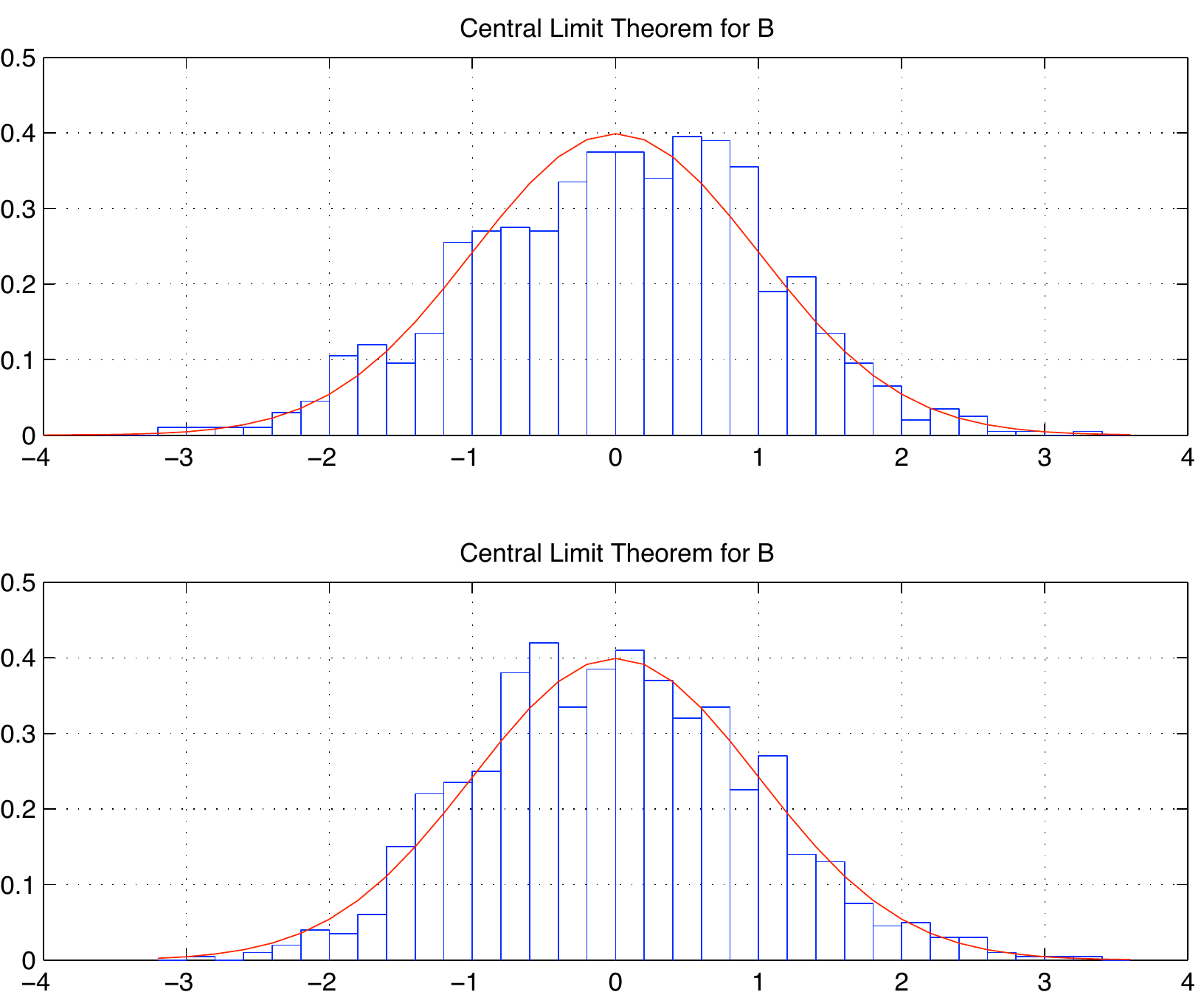}
\caption{Central Limit Theorem}
\end{figure}


\section{CONCLUSION}


Via the use of a persistently excited adaptive tracking control, 
we have shown that it was possible to get ride of the strong controllability 
assumption recently proposed by Bercu and Vazquez \cite{BercuVazq}, \cite{BercuVazquez}.
We have established the almost sure convergence for the LS and WLS estimators in the
multidimensional ARX framework. In addition, we have shown
the residual optimality of the adaptive tracking. Moreover,
both LS and the WLS estimators share the same CLT and LIL. 
We hope that similar analysis could be extended to the
ARMAX framework.


\newpage

\section*{Appendix A}

\begin{center}
{\small PROOF OF LEMMA \ref{MAINLEMMA}}
\end{center}

\renewcommand{\thesection}{\Alph{section}} 
\renewcommand{\theequation}
{\thesection.\arabic{equation}} \setcounter{section}{1}  
\setcounter{equation}{0}


Let $\cA$ and $\cB$ be the
infinite-dimensional diagonal square matrices given by
\begin{equation*}
\cA=\left( 
\begin{array}{ccccc}
\Gamma & 0 & \cdots & \cdots & \cdots \\ 
0 & \Gamma & 0 & \cdots & \cdots \\ 
\cdots & 0 & \Gamma & 0 & \cdots \\ 
\cdots & \cdots & \cdots & \cdots & \cdots \\ 
\cdots & \cdots & \cdots & \cdots & \cdots%
\end{array}
\right),
\end{equation*}

\begin{equation*}
\cB=\left( 
\begin{array}{ccccc}
\Delta & 0 & \cdots & \cdots & \cdots \\ 
0 & \Delta & 0 & \cdots & \cdots \\ 
\cdots & 0 & \Delta & 0 & \cdots \\ 
\cdots & \cdots & \cdots & \cdots & \cdots \\ 
\cdots & \cdots & \cdots & \cdots & \cdots%
\end{array}
\right).
\end{equation*}
Moreover, denote by $\cP$ and $\cQ$ the infinite-dimensional rectangular matrices with $dq$
rows and an infinite number of columns, respectively given, if $p\geq q$, by 
\begin{equation*}
\cP=\left( 
\begin{array}{cccccc}
P_p & P_{p+1} & \cdots & P_{k} & P_{k+1} & \cdots \\ 
P_{p-1} & P_p & \cdots & P_{k-1} & P_{k} & \cdots \\ 
\cdots & \cdots & \cdots & \cdots & \cdots & \cdots \\ 
P_{p-q+2} & P_{p-q+3} & \cdots & P_{k-q+2} & P_{k-q+3} & \cdots \\ 
P_{p-q+1} & P_{p-q+2} & \cdots & P_{k-q+1} & P_{k-q+2} & \cdots%
\end{array}
\right),
\end{equation*}

\begin{equation*}
\cQ=\left( 
\begin{array}{cccccc}
Q_p & Q_{p+1} & \cdots & Q_{k} & Q_{k+1} & \cdots \\ 
Q_{p-1} & Q_p & \cdots & Q_{k-1} & Q_{k} & \cdots \\ 
\cdots & \cdots & \cdots & \cdots & \cdots & \cdots \\ 
Q_{p-q+2} & Q_{p-q+3} & \cdots & Q_{k-q+2} & Q_{k-q+3} & \cdots \\ 
Q_{p-q+1} & Q_{p-q+2} & \cdots & Q_{k-q+1} & Q_{k-q+2} & \cdots%
\end{array}
\right),
\end{equation*}
while, if $p\leq q$, by 
\begin{equation*}
\cP=\left( 
\begin{array}{ccccccc}
P_p & P_{p+1} & \cdots & \cdots & P_{k} & P_{k+1} & \cdots \\ 
\cdots & \cdots & \cdots & \cdots & \cdots & \cdots & \cdots \\ 
P_{1} & P_2 & \cdots & \cdots & P_{k-p+1} & P_{k-p+2} & \cdots \\ 
0 & P_1 & P_2 & \cdots & P_{k-p} & P_{k-p+1} & \cdots \\ 
\cdots & \cdots & \cdots & \cdots & \cdots & \cdots & \cdots \\ 
0 & \cdots & \cdots & 0 & P_1 & P_2 & \cdots%
\end{array}
\right),
\end{equation*}

\begin{equation*}
\cQ=\left( 
\begin{array}{ccccccc}
Q_p & Q_{p+1} & \cdots & \cdots & Q_{k} & Q_{k+1} & \cdots \\ 
\cdots & \cdots & \cdots & \cdots & \cdots & \cdots & \cdots \\ 
Q_{1} & Q_2 & \cdots & \cdots & Q_{k-p+1} & Q_{k-p+2} & \cdots \\ 
Q_0 & Q_1 & Q_2 & \cdots & Q_{k-p} & Q_{k-p+1} & \cdots \\ 
0 & Q_0 & Q_1 & \cdots & Q_{k-p-1} & Q_{k-p} & \cdots \\ 
\cdots & \cdots & \cdots & \cdots & \cdots & \cdots & \cdots \\ 
0 & \cdots & \cdots & 0 & Q_0 & Q_1 & \cdots%
\end{array}
\right).
\end{equation*}
Furthermore, let $\Sigma=\Delta-\Delta(\Gamma+\Delta)^{-1}\Delta$ and denote by
$\cC$ the block diagonal matrix of order $dp$ 
\begin{equation*}
\cC=\left( 
\begin{array}{ccccc}
\Sigma & 0 & \cdots & 0 & 0 \\ 
0 & \Sigma & 0 & \cdots & 0 \\ 
\cdots & \cdots & \cdots & \cdots & \cdots \\ 
0 & \cdots & 0 & \Sigma & 0 \\ 
0 & 0 & \cdots & 0 & \Sigma
\end{array}
\right).
\end{equation*}
One can observe that $\Sigma$ is a positive definite matrix.
Finally, if $p\geq q$, denote by $V$ the matrix with $dq$ rows and $dp$ columns given by
\begin{equation*}
V=\left( 
\begin{array}{ccccccc}
D_0 & D_1 & D_2 & \cdots & \cdots & D_{p-2} & D_{p-1} \\ 
0 & D_0 & D_1 & \cdots & \cdots & D_{p-3} & D_{p-2} \\ 
\cdots & \cdots & \cdots & \cdots & \cdots & \cdots & \cdots \\ 
0 & \cdots & D_0 & D_1 & D_2 & \cdots & D_{p-q+1} \\ 
0 & \cdots & \cdots & D_0 & D_1 & \cdots & D_{p-q}%
\end{array}
\right),
\end{equation*}
while, if $p\leq q$, the upper triangular square matrix of order $dp$ given by
\begin{equation*} 
V=\left( 
\begin{array}{ccccc}
D_0 & D_1 & \cdots & D_{p-2} & D_{p-1} \\ 
0 & D_0 & D_1 & \cdots & D_{p-2} \\ 
\cdots & \cdots & \cdots & \cdots & \cdots \\ 
0 & \cdots & 0 & 0 & D_{0} \\ 
\end{array}
\right).
\end{equation*}
On the one hand, if  $p\geq q$, we can deduce from (\ref{SCHUR}) after 
some straightforward, although rather lengthy, linear algebra calculations that 
\begin{equation}  
\label{SCHURDEC}
S=\cP \cA \cP^t+ \cQ \cB \cQ^t+V \cC V^t.
\end{equation}
We shall focus our attention on the last term in (\ref{SCHURDEC}).
Since the matrix $\cC$ is positive definite, it immediately follows that 
$V \cC V^t$ is also positive definite. Consequently, the Schur complement $S$ is invertible.
On the other hand, if $p \leq q$, we can see from (\ref{SCHUR}) that 
\begin{equation*}  
S=\cP \cA \cP^t+ \cQ \cB \cQ^t+R.
\end{equation*}
where $R$ is the symmetric square matrix of order $dq$
\begin{equation*}
R=\left( 
\begin{array}{cc}
V\cC V^t & \cO \\ 
\cO^t & W
\end{array}
\right)
\end{equation*}
where $\cO$ stands for the zeros matrix of order $dp\times d(q-p)$ and $W$ is 
the block diagonal matrix of order $d(q-p)$ 
\begin{equation*}
W=\left( 
\begin{array}{ccccc}
\Delta & 0 & \cdots & 0 & 0 \\ 
0 & \Delta & 0 & \cdots & 0 \\ 
\cdots & \cdots & \cdots & \cdots & \cdots \\ 
0 & \cdots & 0 & \Delta & 0 \\ 
0 & 0 & \cdots & 0 & \Delta
\end{array}
\right).
\end{equation*}
Taking into account the fact that $V\cC V^t$ and $W$ are both positive definite
matrices, we obtain that $R$ is also positive definite which implies that $S$ is invertible.
Finally, we infer from (\ref{DEFLAMBDA}) that
\begin{equation}  \label{DECDET}
\det(\Lambda)=\det(L)\det(S)=\det(\Gamma+\Delta)^p\det(S).
\end{equation}
Consequently, we deduce from (\ref{DECDET}) that $\Lambda$ is invertible and formula 
(\ref{INVLAMBDA}) can be found in \cite{Horn} page 18, which completes the proof of
Lemma \ref{MAINLEMMA}.$\hfill 
\mathbin{\vbox{\hrule\hbox{\vrule height1ex \kern.5em
\vrule height1ex}\hrule}}$


\section*{Appendix B}

\begin{center}
{\small PROOF OF THEOREM \ref{ASPLS}}
\end{center}

\renewcommand{\thesection}{\Alph{section}} 
\renewcommand{\theequation}
{\thesection.\arabic{equation}} \setcounter{section}{2} 
\setcounter{equation}{0}

In order to prove Theorem \ref{ASPLS}, we shall make use of the same 
approach than Bercu \cite{Bercu2} or Guo and
Chen \cite{Guo1}. First of all, we recall that for all $n\geq 0$ 
\begin{equation}
X_{n+1}-x_{n+1}=\pi _{n}+\varepsilon _{n+1}+\xi_{n+1}.  \label{EQBASIS}
\end{equation}
It follows from (\ref{EQBASIS}) together with the strong law of large
numbers for martingales (see e.g. Corollary 1.3.25 of \cite{Duflo}) 
that $n=\mathcal{O}(s_{n})$ a.s. Moreover, by Theorem 1 of \cite{Bercu2} or Lemma 1
of \cite{Guo1}, we have 
\begin{equation}
\sum_{k=1}^{n}(1-f_{k})\parallel \pi _{k}\parallel ^{2}=\mathcal{O}(\log
s_{n})\hspace{0.5cm}\text{a.s.}  \label{SUMPIF}
\end{equation}
where $f_{n}=\Phi _{n}^{t}S_{n}^{-1}\Phi _{n}$. Hence, if 
$(\varepsilon_{n}) $ has finite conditional moment of order $\alpha >2$, we can show by
the causality assumption on the matrix polynomial $B$
together with (\ref{SUMPIF}) that $\parallel \Phi _{n}\parallel ^{2}=
\mathcal{O}(s_{n}^{\beta })$ a.s. for all $2\alpha ^{-1}<\beta <1$. In
addition, let $g_{n}=\Phi _{n}^{t}S_{n-1}^{-1}\Phi _{n}$ and $\delta _{n}=
\text{tr}(S_{n-1}^{-1}-S_{n}^{-1})$. It is well-known that 
\begin{equation*}
(1-f_{n})(1+g_{n})=1
\end{equation*}
and $(\delta _{n})$ tends to zero a.s. Consequently, as 
$$1+g_{n}\leq 2+\delta _{n}\parallel \Phi _{n}\parallel ^{2},$$ 
we infer from from (\ref{SUMPIF}) that 
\begin{equation}
\sum_{k=1}^{n}\parallel \pi _{k}\parallel ^{2}=o(s_{n}^{\beta }\log s_{n})
\hspace{0.5cm}\text{a.s.}  \label{SPI}
\end{equation}
Therefore, we obtain from (\ref{CT}), (\ref{EQBASIS}) and (\ref{SPI}) that 
\begin{equation}
\sum_{k=1}^{n}\parallel X_{k+1}\parallel ^{2}=o(s_{n}^{\beta }\log s_{n})+
\mathcal{O}(n)\hspace{0.5cm}\text{a.s.}  \label{SX}
\end{equation}
Furthermore, as $B$ is causal, we find from relation (\ref{ARX}) that 
\begin{equation}
U_{n}=B^{-1}(R)A(R)X_{n+1}-B^{-1}(R)\varepsilon _{n+1}  \label{SBU}
\end{equation}
which implies by (\ref{SX}) that 
\begin{equation}
\sum_{k=1}^{n}\parallel U_{k}\parallel ^{2}=o(s_{n}^{\beta }\log s_{n})+
\mathcal{O}(n)\hspace{0.5cm}\text{a.s.}  \label{SU}
\end{equation}
It remains to put together the two contributions (\ref{SX}) and (\ref{SU})
to deduce that $s_{n}=o(s_{n})+\mathcal{O}(n)$ a.s. leading to $s_{n}=
\mathcal{O}(n)$ a.s. Hence, it follows from (\ref{SPI}) that 
\begin{equation}
\sum_{k=1}^{n}\parallel \pi _{k}\parallel ^{2}=o(n)\hspace{0.5cm}\text{a.s.}
\label{SP}
\end{equation}
Consequently, we obtain from (\ref{CT}), (\ref{EQBASIS}), (\ref{SP}) 
and the strong law of large numbers for martingales 
(see e.g. Theorem 4.3.16 of \cite{Duflo}) that
\begin{equation*}
\lim_{n\rightarrow \infty }\frac{1}{n}\sum_{k=1}^{n}X_{k}X_{k}^{t}=\Gamma+\Delta 
\hspace{0.5cm}\text{a.s.}
\end{equation*}
and, for all $1\leq i\leq p-1$, 
\begin{equation*}
\sum_{k=0}^{n}X_{k}X_{k-i}^{t}=o(n)
\hspace{0.5cm}\text{a.s.}
\end{equation*}
 which implies that 
\begin{equation}
\lim_{n\rightarrow \infty }\frac{1}{n}
\sum_{k=1}^{n}X_{k}^{p}(X_{k}^{p})^{t}=L\hspace{0.5cm}\text{a.s.}
\label{CVGX}
\end{equation}
where $L$ is given by (\ref{DEFL}). Furthermore, it follows from (\ref{ARX}),
(\ref{EQBASIS}) and (\ref{SBU}) that for all $n\geq 0$ 
\begin{eqnarray*}
U_{n} &=&B^{-1}(R)A(R)X_{n+1}-B^{-1}(R)\varepsilon _{n+1}, \\
&=&V_{n}+W_{n+1}+Z_{n+1},
\end{eqnarray*}
where
\begin{eqnarray*}
V_{n} &=& B^{-1}(R)A(R)(\pi _{n}+x_{n+1}),\\
W_{n+1} &=& P(R)\varepsilon _{n+1},\\
Z_{n+1} &=& B^{-1}(R)A(R)\xi_{n+1}.
\end{eqnarray*}
Consequently, we deduce from the Cauchy-Schwarz inequality together with (\ref{CT}), 
(\ref{SP}), and the strong law of large numbers for martingales (see e.g. Theorem 4.3.16 of \cite{Duflo}) that for all $1\leq i\leq q$
\begin{equation*}
\lim_{n\rightarrow \infty }\frac{1}{n}\sum_{k=1}^{n}U_{k}U_{k-i+1}^{t}=H_{i}
\hspace{0.5cm}\text{a.s.}
\end{equation*}
which ensures that 
\begin{equation}
\lim_{n\rightarrow \infty }\frac{1}{n}
\sum_{k=1}^{n}U_{k}^{q}(U_{k}^{q})^{t}=H\hspace{0.5cm}\text{a.s.}
\label{CVGU}
\end{equation}
where $H$ is given by (\ref{DEFH}). Via the same lines, we also find that
\begin{equation}
\lim_{n\rightarrow \infty }\frac{1}{n}
\sum_{k=1}^{n}X_{k}^{p}(U_{k-1}^{q})^{t}=K^{t}\hspace{0.5cm}\text{a.s.}
\label{CVGXU}
\end{equation}
Therefore, it follows from the conjunction of (\ref{CVGX}), (\ref{CVGU}) and
(\ref{CVGXU}) that 
\begin{equation}
\lim_{n\rightarrow \infty }\frac{S_{n}}{n}=\Lambda \hspace{0.5cm}\text{a.s.}
\label{CVGFIN}
\end{equation}
where the limiting matrix $\Lambda $ is given by (\ref{DEFLAMBDA}). Thanks to Lemma \ref{MAINLEMMA}, the matrix $\Lambda $ is invertible. This is the key point for the rest of the proof. On
the one hand, it follows from (\ref{CVGFIN}) that $n=\mathcal{O}(\lambda
_{min}(S_{n}))$, $\parallel \Phi _{n}\parallel ^{2}=o(n)$ a.s. which implies
that $f_{n}$ tends to zero a.s. Hence, by (\ref{SUMPIF}), we find that 
\begin{equation}
\sum_{k=1}^{n}\parallel \pi _{k}\parallel ^{2}=\mathcal{O}(\log n)\hspace{
0.5cm}\text{a.s.}  \label{PIFIN}
\end{equation}
On the other hand, we obviously have from (\ref{EQBASIS}) 
\begin{equation}
\parallel C_{n}-\Delta _{n}\parallel =\mathcal{O}\left( \frac{1}{n}
\sum_{k=1}^{n}\parallel \pi _{k-1}\parallel ^{2}\right) \hspace{0.5cm}\text{
a.s.}  \label{COSTPI}
\end{equation}
Consequently, we immediately obtain the tracking residual optimality (\ref{TH12})
from (\ref{PIFIN}) and (\ref{COSTPI}). Furthermore, by a well-known result
of Lai and Wei \cite{Lai} on the LS estimator, we also have 
\begin{equation}
\parallel \widehat{\theta }_{n+1}-\theta \parallel ^{2}=\mathcal{O}\left( 
\frac{\log \lambda _{max}S_{n}}{\lambda _{min}S_{n}}\right) \hspace{0.5cm}
\text{a.s.}  \label{RLS}
\end{equation}
Hence (\ref{TH14}) clearly follows from (\ref{CVGFIN}) and (\ref{RLS}), which completes the proof of Theorem \ref{ASPLS}.
$\hfill 
\mathbin{\vbox{\hrule\hbox{\vrule height1ex \kern.5em
\vrule height1ex}\hrule}}$


\section*{Appendix C}

\begin{center}
{\small PROOF OF THEOREM \ref{ASPWLS}}
\end{center}

\renewcommand{\thesection}{\Alph{section}} 
\renewcommand{\theequation}
{\thesection.\arabic{equation}} \setcounter{section}{3} 
\setcounter{equation}{0}

By Theorem 1 of \cite{Bercu1}, we have 
\begin{equation}  \label{SUMPIFA}
\sum_{n=1}^{\infty} a_{n}(1-f_{n}(a))\parallel \pi_{n} \parallel^{2} < +
\infty \hspace{0.5cm} \text{a.s.}
\end{equation}
where the coefficient $f_{n}(a)=a_{n}\Phi_{n}^{t}S_{n}^{-1}(a)\Phi_{n}$. Then, as 
the weighted sequence $(a_n)$ is given by
$$
a_{n}=\Bigl(\frac{1}{\log s_{n}} \Bigr)^{1+\gamma} 
$$ 
with $\gamma>0$, we clearly have 
$a_{n}^{-1}=\mathcal{O}(s_n)$ a.s. Hence, it follows from (\ref{SUMPIFA})
together with Kronecker's Lemma given e.g. by Lemma 1.3.14 of \cite{Duflo}
that 
\begin{equation}  
\label{SPIA}
\sum_{k=1}^{n}\parallel \pi_{k} \parallel^{2} =o(s_{n}) \hspace{0.5cm} 
\text{a.s.}
\end{equation}
Therefore, we obtain from (\ref{CT}), (\ref{EQBASIS}), (\ref{SPIA}) and the strong law of large numbers for martingales (see e.g. Theorem 4.3.16 of \cite{Duflo}) that
\begin{equation}  \label{SXA}
\sum_{k=1}^n\parallel X_{k+1} \parallel^2=o(s_n)+\mathcal{O}(n) 
\hspace{0.5cm} \text{a.s.}
\end{equation}
In addition, we also deduce from
the causality assumption on the matrix polynomial $B$
that 
\begin{equation}  \label{SUA}
\sum_{k=1}^n\parallel U_{k} \parallel^2=o(s_n)+\mathcal{O}(n) \hspace{0.5cm} 
\text{a.s.}
\end{equation}
Consequently, we immediately infer from (\ref{SXA}) and (\ref{SUA}) that $
s_n=o(s_n)+\mathcal{O}(n)$ so $s_n=\mathcal{O}(n)$ a.s. Hence, (\ref{SPIA})
implies that 
\begin{equation}  \label{SPA}
\sum_{k=1}^n\parallel \pi_k \parallel^2=o(n) \hspace{0.5cm} \text{a.s.}
\end{equation}
Proceeding exactly as in Appendix A, we find from (\ref{SPA}) that 
\begin{equation*}
\lim_{n\rightarrow \infty} \frac{S_{n}}{n} = \Lambda \hspace{0.5cm}\text{a.s.
}
\end{equation*}
Via an Abel transform, it ensures that 
\begin{equation}  \label{CVGFINA}
\lim_{n\rightarrow \infty} (\log n)^{1+\gamma} \frac{S_{n}(a)}{n} = \Lambda 
\hspace{0.5cm}\text{a.s.}
\end{equation}
We obviously have from (\ref{CVGFINA}) that $f_{n}(a)$ tends to zero a.s.
Consequently, we obtain from (\ref{SUMPIFA}) and Kronecker's Lemma that 
\begin{equation}  \label{PIFINA}
\sum_{k=1}^{n}\parallel \pi_{k} \parallel^{2} =o((\log s_{n})^{1+\gamma}) 
\hspace{0.5cm} \text{a.s.}
\end{equation}
Then, (\ref{TH22}) clearly follows from (\ref{COSTPI}) and (\ref{PIFINA}).
Finally, by Theorem 1 of \cite{Bercu1} 
\begin{equation}  \label{RWLS}
\parallel \widehat{\theta}_{n+1} - \theta \parallel^{2} =\mathcal{O} \left( 
\frac{1}{\lambda_{min} S_{n}(a)} \right) \hspace{0.5cm} \text{a.s.}
\end{equation}
Hence, we obtain (\ref{TH23}) from (\ref{CVGFINA}) and (\ref{RWLS}), which
completes the proof of Theorem \ref{ASPWLS}. 
$\hfill\mathbin{\vbox{\hrule\hbox{\vrule height1ex \kern.5em\vrule height1ex}\hrule}}$ 


\section*{Appendix D}

\begin{center}
{\small PROOF OF THEOREM \ref{CLTLIL}}
\end{center}

\renewcommand{\thesection}{\Alph{section}} 
\renewcommand{\theequation}
{\thesection.\arabic{equation}} \setcounter{section}{4} 
\setcounter{equation}{0}

First of all, it follows from (\ref{MOD}) and (\ref{WLS}) that for all $n\geq 1$ 
\begin{equation}  
\label{DECWLS}
\widehat{\theta}_{n}-\theta=S_{n-1}^{-1}(a)M_{n}(a)
\end{equation}
where 
\begin{equation}  \label{DEFMART}
M_{n}(a)=\widehat{\theta}_0-\theta+\sum_{k=1}^n
a_{k-1}\Phi_{k-1}\varepsilon_k^t.
\end{equation}
We now make use of the CLT for multivariate martingales given e.g. by Lemma
C.1 of \cite{Bercu2}, see also \cite{Duflo}. On the one hand,
for the LS algorithm, we clearly deduce (\ref{TH31}) from convergence 
(\ref{TH11}) and decomposition (\ref{DECWLS}). On the other hand, for the WLS
algorithm, we also infer (\ref{TH31}) from convergence (\ref{TH21}) and 
(\ref{DECWLS}). Next, we make use of the LIL for multivariate martingales given
e.g. by Lemma C.2 of \cite{Bercu2}, see also \cite{Duflo}, \cite{Stout}. For
the LS algorithm, since $(\varepsilon_n)$ has finite conditional moment of
order $\alpha>2$, we obtain from Chow's Lemma given e.g. by Corollary 2.8.5
of \cite{Stout} that for all $2< \beta < \alpha$ 
\begin{equation}  \label{NORMEPS}
\sum_{k=1}^{n}\parallel \varepsilon_{k}\parallel^{\beta}=\mathcal{O}(n) 
\hspace{0.5cm}\text{a.s.}
\end{equation}
The exogenous noise $(\xi_n)$ shares the same regularity in norm than
$(\varepsilon_n)$ which means that for all $2< \beta <\alpha$ 
\begin{equation}  
\label{NORMXI}
\sum_{k=1}^{n}\parallel \xi_{k}\parallel^{\beta}=\mathcal{O}(n) 
\hspace{0.5cm}\text{a.s.}
\end{equation}
Consequently, as the reference trajectory $(x_n)$ satisfies (\ref{REGNORM}),
we deduce from (\ref{EQBASIS}) together with (\ref{PIFIN}), (\ref{NORMEPS}) and
(\ref{NORMXI})
that for some $2< \beta < \alpha$
\begin{equation}  \label{NORMX}
\sum_{k=1}^{n}\parallel X_{k}\parallel^{\beta}=\mathcal{O}(n) \hspace{0.5cm} 
\text{a.s.}
\end{equation}
Furthermore, it follows from (\ref{SBU}) and (\ref{NORMX}) that 
\begin{equation}  \label{NORMU}
\sum_{k=1}^{n}\parallel U_{k}\parallel^{\beta}=\mathcal{O}(n) \hspace{0.5cm} 
\text{a.s.}
\end{equation}
Hence, we clearly obtain from (\ref{NORMX}) and (\ref{NORMU}) that 
\begin{equation}  \label{NORMPHI}
\sum_{k=1}^{n}\parallel \Phi_{k}\parallel^{\beta}=\mathcal{O}(n) \hspace{
0.5cm}\text{a.s.}
\end{equation}
Therefore, as $\beta >2$, (\ref{NORMPHI}) immediately implies that 
\begin{equation*}
\sum_{n=1}^{\infty}\left(\frac{\parallel \Phi_{n}\parallel} {\sqrt{n}}
\right)^{\beta}<+\infty \hspace{0.5cm}\text{a.s.}
\end{equation*}
Finally, Lemma C.2 of \cite{Bercu2} together with convergence (\ref{TH11})
and (\ref{DECWLS}) lead to (\ref{TH32}). The proof for the WLS algorithm is
left to the reader because it follows essentially the same arguments than
the proof for the LS algorithm. It is only necessary to add the weighted
sequence $(a_n)$ and to make use of convergence (\ref{TH21}).
$\hfill\mathbin{\vbox{\hrule\hbox{\vrule height1ex \kern.5em\vrule height1ex}\hrule}}$ 


\end{document}